\documentclass[11pt,twoside]{article}

\usepackage{mathrsfs,amsfonts,amsmath,amssymb}
\usepackage{latexsym}
\usepackage{graphicx}
\newtheorem{theorem}{Theorem}

\renewcommand{\thefootnote}{}

\def\R{\mathbb {R}}
\def\E{\mathsf{E}}

\def\({\big (}
\def\){\big )}

\def\_phi{\varphi}

\def\E{{\sf E}}

\author{M. Rudnev\footnote{
		Partially supported by the Leverhulme  Trust Grant RPG--2017--371.} }
\title{ On the number of hinges defined by a point set in $\mathbb R^2$}
	
\date{}
\begin{document}
	\maketitle

\begin{abstract}
It is shown that the number of distinct types of three-point hinges, defined by a real plane set of $n$ points is $\gg n^2\log^{-3} n$, where a hinge is identified by fixing two pair-wise distances in a point triple.
This is achieved via strengthening (modulo a $\log n$ factor) of the Guth-Katz estimate for the number of pair-wise intersections of lines in $\R^3$, arising in the context of the plane Erd\H os distinct distance problem, to a second moment incidence estimate. This relies, in particular, on the generalisation of the Guth-Katz  incidence bound by Solomon and Sharir.  \end{abstract}

\renewcommand{\thefootnote}{\arabic{footnote}}
\section{Introduction}
Given an $n$-point set $\mathcal P\subset \R^2$, define a hinge as an equivalence class $h=[p,q,r]$ of $(p,q,r)\in \mathcal P^3$, identified by a pair of two fixed distances $\|p-q\|,\,\|q-r\|$. Let $H=H(\mathcal P)$ be the set of hinges and $r_H(h)$ the number of realisations of a hinge $h\in H(\mathcal P)$, that is the number of triples $(p,q,r)$ in the equivalence class $h$. What is the minimum number of hinges defined by a $n$-point set? The example of the square grid suggests that it should be $\geq c \,n^2\log^{-1} n$, for some universal $c>0$. We show here that this is indeed the right asymptotics, up to an additional factor of $\log^2 n$.

The question on what is the minimum number of hinges defined by a point set was asked, to the author's knowledge,  by Iosevich and Passant \cite{IP}. In general, what other conclusions can be drawn apropos of distances in the plane in the wake of the breakthrough resolution of the Erd\H os distinct distance conjecture \cite{Erd} by Guth and Katz \cite{GuK}, within the scope of techniques, founded in the latter paper? It can be easily inferred from the main incidence estimate of \cite{GuK} (estimate \eqref{pkb} below) that the number of distinct congruence classes of {\em triangles}, defined by an essentially non-collinear plane $n$-point set is $\geq c\,n^2$ \cite{Rud}.

A distance is an equivalence class of pairs of points in the set. Its natural generalisation is a {\em chain}, that is a $k$-tuple of points with $k-1\geq 1$ fixed distances along the tuple. For $k>1$ it is not a rigid configuration, which accounts for higher volatility of the number of realisations function for chains versus rigid point configurations. It was observed in \cite{PSS} that the maximum number of realisations of a hinge, that is a two-chain, can be $\geq c\, n^2$ -- just take about half of the points on each of two concentric circles and add the centre. This is much greater than what is believed to be the case with rigid point configurations, where it is conjectured that the maximum number of realisations of a single distance is at most $n^{1+\epsilon}$, for any $\epsilon >0$ and a sufficiently large $n$.

Naturally, questions about realisation functions of point configurations with prescribed metric structure have continuous and finite field analogues -- see e.g. \cite{IoP} for some discussion and references.

\medskip
\noindent
To proceed with our consideration of hinges, define 
$$\E_H(\mathcal P):= \sum_{h\in H} r^2_H(h).$$
Two triples $(p,q,r)$ and $(p',q',r')$ are in the same equivalence class if and only if simultaneously
\begin{equation}\label{hc}
\|p-q\|=\|p'-q'\|,\qquad \|r-q\|=\|r'-q'\|.
\end{equation}
Consider the Elekes-Sharir map (see \cite{ES}, \cite{GuK} and  \cite{RS} for some generalisations) $\phi:\,\mathcal P^2\to \mathcal K(\R^3)$, acting as $(p,p')\to l_{pp'}$, where $\mathcal K(\R^3)$ is the Klein quadric, the set of lines in $\R \mathbb P^3$. Explicitly, in Pl\"ucker coordinates (see \cite{RS} or \cite{ES}, \cite{GuK} for the expressions as actual lines in $\R^3$) with $p=(p_1,p_2)$, etc., one has
\begin{equation}
l_{pp'}=\left(\frac{p_2'-p_2}{2}:\frac{p_1-p_1'}{2}:1:\frac{p_2'+p_2}{2}:-\frac{p_1+p_1'}{2}:\frac{\|p\|^2-\|p'\|^2}{4}\right)\,.
\label{lines}\end{equation}
It is well known (\cite{GuK}, and  \cite{RS} for some generalisations)  that  the set of $n^2$ lines $L:=\{l_{pp'}\}_{(p,p')\in \mathcal P^2}$ has the property that just $O(n)$ may be concurrent, coplanar or lie in a regulus. 

The hinge condition \eqref{hc} holds if and only if simultaneously 
$$
l_{pp'}\cap l_{qq'}\neq \emptyset, \qquad l_{rr'}\cap l_{qq'}\neq \emptyset.
$$ 
It follows that, $\nu(l)$ denoting the number of other lines in $L$, meeting some $l\in L$, 
$$
\E_H(\mathcal P)\ll \sum_{l\in L} \nu^2(l)\,.
$$

This note uses the standard $\ll,\gg,\sim$ notations, as well as $O(\cdot)$ for $\ll$, to subsume absolute constants. All point/line sets involved  are finite, of cardinality $|\cdot|$.

\medskip
Its main (and only) claim is as follows.

\begin{theorem} One has $\E_H(\mathcal P)\ll n^4 \log^3n$, hence $|H|\gg n^2\log^{-3}n$, for any $\mathcal P\subset \R^2$, with $|\mathcal P|=n,$ sufficiently large.
\label{hinge}\end{theorem}

The proof in the next section consists in multiple applications of the Guth-Katz type incidence bounds for lines in $\R^3$. These include the following Guth-Katz vintage incidence bound\footnote{The following formulation, see e.g. \cite[Theorem 1.1]{SS} is not given explicitly in \cite{GuK} but follows easily.}, as well as some of its more recent developments.

\begin{theorem}[Guth-Katz] \label{GK}  Let $P$ be a set of points and $L$ a set of lines in $\R^3$,
let $ s \leq |L|$  be a parameter, such that no affine plane contains more than $s$
lines of $L$. The number $I(P,L)$ of incidences between $P$ and $L$ satisfies the bound
$$I(P,L) \ll |L|^{3/4}|P|^{1/2}+  s^{1/3} |L|^{1/3}|P|^{2/3} +|L|+|P|\,.$$
\end{theorem}
The considerations in \cite{GuK} imply that a typical line from $n^2$ lines $l_{pp'}\in L$ from \eqref{lines} meets $O(n \log n)$ other lines.
Theorem \ref{hinge} extends this to a second moment bound: the average, over lines, of the square of the number of lines a line meets, is $O(n^2 \log^3 n)$, for the extra price of a $\log n$ factor.

 By taking $\mathcal P$ to be a square grid, it's easy to see from \eqref{lines} that a typical one of the $n^2$ lines in $L$ meets $\sim n\log n$ lines\footnote{Two lines $(a:b:1:c:-d:-ac+bd)$ and  $(a':b':1:c':-d':-a'c'+b'd')$ in Pl\"ucker coordinates meet iff $ac'-bd'+a'c-b'd = ac-bd+a'c'-b'd'$. In the example in question $a,\ldots,d'$ are half-integers $O(\sqrt{n})$, and after rearranging and changing variables one sees that for a typical $(a,b,c,d)$, the number of quadruples $(a',b',c',d')$ satisfying the latter equation is roughly the number of quadruples of natural numbers $n_1/n_2=n_3/n_4$, with $n_i=O(\sqrt{n})$, which is known and easily seen to be $\sim n\log n$.}, and no line meets more than $O(n\log n)$ other lines. So the estimate of Theorem \ref{hinge} may be off by $\log n$ from being sharp, while the incidence Theorem \ref{GK} is sharp.

\medskip
To establish Theorem \ref{hinge}, one uses two somewhat more elaborate variants of Theorem \ref{GK}, as follows. The first one is a bipartite version of \cite[Theorem 4.1]{GuK} (with the two line sets involved not being necessarily disjoint). A proof can be  extracted from the author's paper \cite[Proof of Theorem 12]{Ru}. The latter was followed by de Zeeuw's work \cite{FdZ}, which offers an almost\footnote{De Zeeuw uses the term {\em quadric} for both plane and regulus. Furthermore, his quadric condition, in comparison to the formulation here  is less strict: no quadric may contain $s$ points of {\em each} line family, hence instead of the term $s|L'|$ here he has $s(|L|+|L'|)$. The fact that the latter term gets replaced by $s|L'|$ is immediate from  \cite[Proof of Lemma 3.1, last passage on p.3]{FdZ}.} verbatim version of the statement below  \cite[Lemma 3.1]{FdZ}, and proof.

\begin{theorem} \label{GKb} Let $L,L'$ be two sets of lines in $\R^3$, with $|L'|\leq |L|$ and at most $s\leq |L|$ of lines from $L$ lying in a plane or regulus. If $P$ is a set of all points where two distinct lines $l,l'$, with $l\in L$, $l'\in L'$ meet, then
$$|P|\ll |L|\sqrt{|L'|} + s|L'|\,.$$
\end{theorem}

Theorem \ref{GK} was generalised by Sharir and Solomon \cite[Theorem 1.3(a)]{SS} as follows.
\begin{theorem} \label{GKi} Let $P$ be a set of points and $L$ a set of lines in $\R^3$, lying in a degree $d\geq 2$ polynomial surface, which does not contain linear factors. Suppose, at most $s\leq d$ lines\footnote{The restriction  $s\leq d$ is merely a consequence of the Bezout theorem; certainly $s$ can be replaced by a bigger quantity. Also note that any set  of lines $L$ can be included into a zero set of a polynomial of degree $d\ll \sqrt{|L|}$. If $d\sim \sqrt{|L|}$, the bounds of Theorems \ref{GK} and \ref{GKi} coincide and absorb the extra term $|L|d$ that would come from dropping the assumption on the lack of linear factors in the latter theorem. } are contained in any plane. The number $I(P,L)$ of incidences between $P$ and $L$ satisfies the bound
\begin{equation}\label{ssbound}I(P,L) \ll |P|^{1/2}|L|^{1/2}d^{1/2} + |P|^{2/3}d^{2/3}s^{1/3} + |P| +|L|\,.\end{equation}
\end{theorem}
Technically, this is largely Theorem \ref{GKi} that allows for the result here, for its proof in some sense iterates the Guth-Katz argument, applying the incidence bound of Theorem \ref{GK} only to a $O(d^2)$-strong subset of lines \cite[Theorem 3.4]{SS}. 

\section{Proof of Theorem \ref{hinge}}

\begin{proof}
Partition, for dyadic $k=2,4,8,\ldots,\leq n$ the set $P$ of pair-wise intersections of lines in the set of lines $L$, defined by \eqref{lines} into sets $P_k$, consisting of intersection points where some number of lines in the interval $ [k,2k)$ lines meet. 

The key result of \cite{GuK} is the following bound (with $|L|=n^2$, $k,s\ll n$ -- this bound combines the estimates of Theorems \ref{GK} and \ref{GKb} with $L=L'$): \begin{equation}\label{pkb} |P_k|\ll \frac{n^3}{k^2}\,.\end{equation}

Given $k$, for dyadic $t=1,2,4,\ldots,\ll \frac{n^2}{k}$, let $L_{k,t}$ be the set of $t$-rich lines, relative to points $P_k$, that is the set of lines in $L$, supporting some number of points in $P_k$ in the interval $ [t,2t)$.

Thus each line $l\in L$ can belong to $O(\log n)$ sets $L_{k,t}$. If  $l$ supports $\nu(l)$ line intersections altogether, then partitioning these intersections by sets $P_k$, with $t_k$ intersection points lying in $P_k$, yields
$$
\nu^2(l)\ll \log n \sum_k (k t_k)^2.
$$
To proceed with the proof of Theorem \ref{hinge}, we assume, for contradiction, can be violated only if there is some pair $(k,t)$, such that
\begin{equation} \label{contr}
|L_{k,t}|(kt)^2 \geq C n^4 \,,
\end{equation}
for a sufficiently large $C$.

Since trivially $|L_{k,t}|\leq n^2$, \eqref{contr} may possibly take place only for $kt \geq n\sqrt{C}$, in particular $t>1$ being sufficiently large.

\medskip
The rest of the proof is a case-by-case analysis, with \eqref{contr} leading to a contradiction in each case.

First, consider separately the case $k=O(1)$ by applying Theorem \ref{GKb} to $L$ and $L'=L_{k,t}$.
It follows that
$$
t |L_{k,t}| \ll n^2 \sqrt{|L_{k,t}| }+ n|L_{k,t}|\,.
$$
If the second term in the latter inequality dominates, then $t=O(n)$, and since $k=O(1)$, \eqref{contr} cannot hold. If the first term dominates, then
$$
|L_{k,t}| \ll \frac{n^4}{t^2}\,,$$
and once again  \eqref{contr} gets contradicted. 

\medskip
From now on we may assume that $k$ is sufficiently large, relative to absolute constants, hidden in the $\ll,\gg$ symbols.

Restrict $P_k$ to only its points supported on lines in $L_{k,t}$. Apply Theorem \ref{GK} to estimate the number of incidences $I(P_k,L_{k,t})$:
\begin{equation}
t |L_{k,t}| \leq  I(P_k,L_{k,t}) \ll |P_k|^{1/2}|L_{k,t}|^{3/4} + n^{1/3}|P_k|^{2/3}|L_{k,t}|^{1/3} + |P_k|.
\label{inc}\end{equation}
There are three cases to consider. 

\medskip {\sf Case (i).}  Suppose $t |L_{k,t}| \ll |P_k|^{1/2}|L_{k,t}|^{3/4}$. It follows that
$$
|L_{k,t}|(kt)^2 \ll \frac{|P_k|^2}{t^4}(kt)^2\ll \frac{n^6}{(kt)^2} \leq n^4\,,
$$
using \eqref{pkb} and that implicit in \eqref{contr}, $kt\geq n$. Thus  \eqref{contr} does not hold in Case (i).

\medskip {\sf Case (ii).}
Now suppose the third term in \eqref{inc} dominates, namely
\begin{equation}\label{case2}
|L_{k,t}|\ll \frac{|P_k|}{t}\,.
\end{equation}
Hence, using \eqref{pkb},
$$
|L_{k,t}|(kt)^2\ll {n^3}t\,,
$$ therefore \eqref{contr} may possibly be true for $t\gg Cn$ only.

We proceed by putting the set $L_{k,t}$ in a zero set $Z$ of a polynomial of degree 
$$d\ll \min\left\{\sqrt{ \frac{|P_k|}{t}},\,n\right\}\,,$$
so $P_k\subset Z$, considering incidences between the set of lines $L$ and $P_k$ (recall that $P_k$ has been restricted to points lying on lines in $L_{k,t}$). Let us partition $L$ into $L^\perp$ and $L^{\|}$, where members of $L^\perp$ do not lie in $Z$, while those of $L^{\|}$ do. 


Since a line in  $L^\perp$ cannot meet $Z$ at more than $d$ points,
\begin{equation}
k|P_k| \leq I(P_k,L)\ll  n^2 d + I(P_k, L^{\|})\,.
\label{bez}\end{equation}
If the first term in the estimate dominates, then by the estimate on $d$ and \eqref{case2} one gets
$$
t\ll \frac{n^4}{k^2|P_k|} \ll \frac{n^4}{k^2 t |L_{k,t}|}\,,
$$ and therefore \eqref{contr} is contradicted.

Otherwise assume (adding, for convenience,  $|P_k|$ on the next line is harmless by \eqref{bez})
\begin{equation}\label{ass}
n^2 d +|P_k| \ll I(P_k, L^{\|})\,.
\end{equation}

Next, we can remove from $Z$ the union of linear factors, because the total contribution of incidences, supported on lines therein is $\ll n^2d+|P_k|$. This follows readily from the trivial estimate $\ll n^2$ for the number of incidences inside each plane factor (since it must have $\ll n$ lines) and the fact that there are at most $d$ of these factors.

Proceeding under the assumption that $Z$ has no plane factors, estimate \eqref{ssbound}, since $k$ is sufficiently large, yields
\begin{equation}\label{ssapp}
I(P_k, L^{\|}) \ll |P_k|^{1/2}d^{1/2}n + |P_k|^{2/3}d^{2/3}n^{1/3} + n^2\,.
\end{equation}
If the first term in the latter estimate dominates, then by \eqref{ass} one has $d \ll |P_k|/n^2$, and hence, from \eqref{bez}, $k|P_k| \ll |P_k|$. This is a contradiction, for $k$ is meant to be sufficiently large. 

If the third term dominates, then $|P_k|\ll \frac{n^2}{k}$ and from \eqref{case2}
\begin{equation}\label{tr}|L_{k,t}| (kt)^2 \ll n^2 (kt) \leq n^4\,,\end{equation}
 since clearly $kt\leq n^2$.

Thus  it remains to consider the dominance of the second term in  estimate \eqref{ssapp}.
If so, we would have 
$$
k|P_k|,\, n^2 d  \ll |P_k|^{2/3}d^{2/3}n^{1/3} \,.
$$
From the second inequality 
$$
d^{1/3} \ll n^{-5/3} |P_k|^{2/3}\,,
$$
so
$$
k|P_k| \ll n^{-3} |P_k|^2,
$$
thus $|P_k|\gg {n^3}{k},$ which contradicts \eqref{pkb}. Thus we conclude that \eqref{contr} does not hold in Case (ii).

\medskip
{\sf Case (iii).}  The term $n^{1/3}|P_k|^{2/3}|L_{k,t}|^{1/3} $ dominates estimate \eqref{inc}, that is 
\begin{equation}
|L_{k,t}|\ll |P_k|n^{1/2}t^{-3/2} \ll \frac{n^{7/2}}{k^2 t^{3/2}}\,.
\label{case3}\end{equation}
Thus \eqref{contr} can possibly hold only if $t\gg C^2n$.

We repeat the analysis from Case (ii), now with $d\ll |P_k|^{1/2}n^{1/4}t^{-3/4}.$ Returning to \eqref{bez}, if the first term in the right-hand side dominates we have

$$
k|P_k|\ll n^2d\ll |P_k|^{1/2}n^{9/4}t^{-3/4}\,.
$$
Hence 
$$
|P_k|^{2/3} \ll n^3 t^{-1}k^{-4/3}, \qquad t|L_{k,t}|^{2/3} \ll n^{10/3}t^{-1} k^{-4/3}\,,
$$
the latter by \eqref{inc}. Thus
$$
|L_{k,t}|(kt)^2 \ll \frac{n^5}{t}\leq n^4\,,
$$
given that $t\geq n$.

The rest of the analysis repeats what has already been done in Case (ii) apropos of the relations \eqref{bez}--\eqref{ssapp}. The only small deviation from the above argument is that  dominance of the term $n^2$ in estimate \eqref{ssapp} would mean, once again, $|P_k|\ll \frac{n^2}{k}$, and thus by \eqref{case3}, instead of estimate \eqref{tr} one has
$$
|L_{k,t}|(kt)^2 \ll n^{5/2}(kt)t^{-1/2} \leq n^4\,,
$$ since $t\geq n$ and $kt\leq n^2$.

Thus assuming \eqref{contr} results in a contradiction in Case (iii) as well, which concludes the proof of Theorem \ref{hinge}. $\hfill \Box$
\end{proof}

\vspace{1cm}
\noindent {\sf Misha Rudnev}\\
School of Mathematics, Fry Building, Bristol BS8 1TH, UK.\\
{\tt misharudnev@gmail.com}

\end{document}